\documentclass[12pt,draft]{amsart}
\usepackage{amssymb, amsmath}
\usepackage{url}
\newtheorem{theorem}{Theorem}[section]

\newtheorem{prop}[theorem]{Proposition}

\theoremstyle{remark}
\newtheorem{remark}[theorem]{Remark}

\theoremstyle{definition}
\newtheorem{defn}[theorem]{Definition}
\newtheorem{example}[theorem]{Example}

\numberwithin{equation}{section}
\numberwithin{theorem}{section}

\newcommand{\intab}{\int_a^b}
\newcommand{\Bp}{{\mathcal B}_p(\R)}
\newcommand{\Ap}{{\mathcal A}_p(\R)}
\newcommand{\Psif}{\Psi_{\!f}}
\newcommand{\est}{e^{-ist}}
\newcommand{\FT}{\mathcal{F}}

\let\abs=\envert

\newcommand{\gonehat}{\widehat{g_1}}
\newcommand{\gtwohat}{\widehat{g_2}}

\newcommand{\hatT}{\hat{T}}

\newcommand{\fhat}{\hat f}
\newcommand{\ghat}{\hat g}

\newcommand{\psihat}{\hat\psi}
\newcommand{\phihat}{\hat \phi}
\newcommand{\Lone}{\Lany{1}}

\newcommand{\intinf}{\int^\infty_{-\infty}}

\newcommand{\N}{{\mathbb N}}
\newcommand{\R}{{\mathbb R}}

\newcommand{\Sc}{{\mathcal S}(\R)}
\newcommand{\Scp}{{\mathcal S}'(\R)}
\newcommand{\fn}{\!:\!}

\newcommand{\llim}{\lim\limits}

\newcommand{\bv}{{\mathcal BV}}

\providecommand{\abs}[1]{\lvert#1\rvert}
\providecommand{\norm}[1]{\lVert#1\rVert}
\providecommand{\Lany}[1]{L^{#1}(\R)}

\begin{document}


\date{Preprint February 25, 2025.  To appear in {\it Czechoslovak Mathematical Journal}}

\title[fourier transform]{the fourier transform in lebesgue spaces}
\author{Erik Talvila}
\address{Department of Mathematics \& Statistics\\
University of the Fraser Valley\\
Abbotsford, BC Canada V2S 7M8}
\email{Erik.Talvila@ufv.ca}

\begin{abstract}
For each $f\in L^p({\mathbb R)}$ ($1\leq p<\infty$) 
it is shown that 
the Fourier transform is the distributional derivative
of a H\"older continuous function.  For each $p$ a norm is defined
so that the space Fourier transforms is isometrically isomorphic
to $L^p({\mathbb R)}$.  There is an exchange theorem and inversion
in norm.
\end{abstract}

\dedicatory{Dedicated to the memory of Jaroslav Kurzweil.}
\maketitle

\section{Introduction}\label{sectionintroduction}
In this paper we use differentiation and integration in spaces of tempered 
distributions to define the Fourier transform of $f\in\Lany{p}$ for all
$1\leq p<\infty$ as the distributional derivative of a continuous function.

If $f\fn\R\to\R$ its Fourier transform is $\FT[f](s)=\fhat(s)=\intinf \est f(t)\,dt$
for $s\in\R$.  The Fourier transform is a bounded linear operator, $\FT\fn\Lone\to\Lany{\infty}$
and $\norm{\fhat}_\infty\leq\norm{f}_1$.  When $f\in\Lone$ the Fourier transform is known to be a continuous
function that vanishes at infinity (Riemann--Lebesgue lemma).  When $f\in\Lany{p}$ for
some $1<p\leq 2$ then interpolation methods can be used to define the Fourier transform.
In this case, the Hausdorff--Young--Babenko--Beckner inequality states that
$\norm{\fhat}_q\leq B_q\norm{f}_p$.
Here, $q$ is the conjugate exponent, $1/p+1/q=1$.  
The constant $B_q$ is given in Remark~\ref{remarkBabenko} below.
The Fourier transform is then a bounded linear
operator from $\Lany{p}$ to $\Lany{q}$.  

For the necessary Fourier transform background
see, for example, \cite{folland}, \cite{grafakos}, \cite{liebloss} and \cite{steinweiss}.

Interpolation methods are not applicable
when $p>2$.  Indeed, see \cite{abdelhakim} where examples are given to show the Fourier transform
is unbounded as an operator to any $\Lany{q}$ space.

For each $f\in\Lany{p}$ ($1\leq p\leq \infty$) the Fourier transform can be defined
as a distribution via the exchange formula $\langle \fhat,\phi\rangle=\langle f,\phihat\rangle
=\intinf f(s)\phihat(s)\,ds$.  Here, $\phi\in\Sc$ where the Schwartz space consists of those
functions in $C^\infty(\R)$ such that for all integers
$m,n\geq 0$,
$x^m\phi^{(n)}(x)\to 0$ as $\abs{x}\to\infty$.  When $p>2$ it may happen that $\fhat$ is not a function.
See, for example, \cite[4.13, p.~34]{steinweiss}.

For background on distributions see, for example, \cite{donoghue}, \cite{folland}, 
\cite{friedlanderjoshi}, \cite{grafakos},
\cite{steinweiss}.

If $f\in\Lone$ then we can integrate the Fourier transform to define 
$\Psif(s)=\int_0^s\fhat(\sigma)\,d\sigma=
\intinf\left(\frac{1-e^{-ist}}{it}\right) f(t)\,dt$.  Interchange of orders of integration
follows by the Fubini--Tonelli theorem.  In this case, the pointwise derivative
$\Psif'(s)=\fhat(s)$ at each point $s\in\R$ and the derivative can be computed using
the fundamental theorem of calculus since $\fhat$ is continuous.  As well, differentiation
under the integral sign is valid by the dominated convergence theorem and yields the same
result.  When $f\in\Lany{p}$ for any $1<p<\infty$ this differentiation formula is
used as the definition of the Fourier transform, the derivative now existing in the
distributional sense.

The kernel $\Psif$ appears several places in the literature connected with Fourier transforms.
It has been used by Titchmarsh \cite{titchmarshcontribution} and Bochner 
\cite{bochner}, although only
for $1\leq p\leq2$ and only for pointwise derivatives.  
For an application to distributions,
see \cite[\S31]{donoghue}.

The outline of the paper is as follows.
In Theorem~\ref{theoremPsif}
the inequality
$\abs{\Psif(s)}\leq C_q\norm{f}_p\abs{s}^{1/p}$ is proved and it is shown
that $\Psif$ is H\"older continuous.  Two Banach spaces are defined in 
Theorem~\ref{theoremBanachspaces}, both isometrically isomorphic to $\Lany{p}$.
The first consists of the functions $\Psif$ for all $f\in\Lany{p}$ and the
second is the set of Fourier transforms of $\Lany{p}$ functions.  Various
properties of the Fourier transform are proved in Theorem~\ref{theoremalgebraic}.
The Fourier transform of a function in $\Lany{p}$ is a tempered distribution.
Definitions of local and global integration are given in 
Definition~\ref{defnintegrationinAp}.  Section~\ref{sectionexchange}
contains an exchange theorem and inversion using a summability kernel.
Fourier transforms of convolutions are considered in Section~\ref{sectionconvolution}.

\section{Properties of $\Psif$}
Various properties of the kernel $\Psif$ are proved here.  A sharp pointwise estimate is 
given and it is shown that
$\Psif$ is H\"older continuous and need not be of locally bounded variation.

\begin{theorem}\label{theoremPsif}
Let $f\in\Lany{p}$ for some $1\leq p<\infty$.  Define $\Psif(s)=
\intinf\left(\frac{1-e^{-ist}}{it}\right) f(t)\,dt$ for $s\in\R$.
Let $1/p+1/q=1$.  If $p=1$ then $q=\infty$.
Let $C_q=4^{1/q}(\int_0^\infty\abs{\sin(t)/t}^q\,dt)^{1/q}$,
with $C_\infty=1$.

(a)  For all $s\in\R$,
$\abs{\Psif(s)}\leq C_q\norm{f}_p\abs{s}^{1/p}$.

(b) For $1<p<\infty$, this estimate is sharp in the limit $\abs{s}\to\infty$ in the sense
that if $\psi\fn\R\to(0,\infty)$ is any function such that $\psi(s)=o(\abs{s}^{1/p})$
as $\abs{s}\to\infty$ then there is $f\in\Lany{p}$ such that
$\Psif(s)/\psi(s)$ is not bounded as $\abs{s}\to\infty$.
In a similar manner the estimate is sharp in the limit $s\to0$.
Also, the constant $C_q$ cannot be replaced with any smaller number.

(c) For $p=1$ the estimate $\abs{\Psif(s)}\leq \norm{f}_p\abs{s}$ is sharp
in the limit $s\to0$.  And, $\Psif(s)=o(s)$ as $\abs{s}\to\infty$.

(d) $\Psif$ is H\"older continuous on $\R$ with exponent $1/p$
and hence uniformly continuous on $\R$.  When $p=1$ the function $\Psif$ is
Lipschitz continuous on $\R$.
\end{theorem}
\begin{proof}
(a) Define $u_s(t)=\frac{1-e^{-ist}}{it}$ with $u_s(0)=s$.  A Taylor expansion of the
exponential function shows $u_s$ is real analytic on $\R$, for each fixed
$s\in\R$.  Note that
\begin{equation}
\abs{u_s(t)}^2=\frac{(1-\cos(st))^2+\sin^2(st)}{t^2}=\frac{2(1-\cos(st))}{t^2}=
\frac{4\sin^2(st/2)}{t^2}.\label{u_s}
\end{equation}
The H\"older inequality now shows $\abs{\Psif(s)}\leq C_q\norm{f}_p\abs{s}^{1/p}$,
with $C_\infty=1$.

(b) To prove the estimate is sharp, define the family of linear operators
$L_s\fn\Lany{p}\to\R$ by $L_s[f]=\Psif(s)/(C_q\psi(s))$ for each $s\in\R$.
For each $s\in\R$ and each $f\neq 0$ the estimate $\abs{L_s[f]}/\norm{f}_p\leq \abs{s}^{1/p}/\psi(s)$
shows $L_s$ is a bounded linear operator.  Now show this family is not uniformly
bounded.  We can write $u_s(t)=\sin(st)/t-i(1-\cos(st))/t=\abs{u_s(t)}e^{i\theta(st)}$
for some phase $-\pi<\theta\leq\pi$ that is a function of the product $st$.
Define $f_s\in\Lany{p}$ by $f_s(t)=\abs{u_s(t)}^{q/p}e^{-i\theta(st)}$.  By the condition
for equality in the H\"older inequality,
$\abs{L_s[f_s]}/\norm{f_s}_p=\abs{s}^{1/p}/\psi(s)$.  This is not bounded as $\abs{s}\to\infty$.
Hence, the family $\{L_s\mid s\in\R\}$ is not uniformly bounded.  By the
Uniform Boundedness Principle there is $f\in\Lany{p}$ such that
$L_s[f]=\Psif(s)/(C_q\psi(s))$ is not bounded as $\abs{s}\to\infty$.
The inequality is then sharp in this limit.  The case $s\to 0$ is similar.  The example
$f_s$ shows $C_q$ cannot be replaced with any smaller number.

(c) Let $f=\chi_{(-1,1)}$.  Then
$\Psif(s)=2\int_0^s(\sin(\sigma)/\sigma)\,d\sigma$.
As $s\to0$, we have $\Psif(s)\sim 2s=\norm{f}_1s$.  The estimate is sharp in
the limit $s\to0$.

Now let $f\in\Lone$.  Let $0<M<s$.  Then $\Psif(s)=\int_0^M\fhat(\sigma)\,d\sigma+
\int_M^s\fhat(\sigma)\,d\sigma$.  Note that $\abs{\int_0^M\fhat(\sigma)\,d\sigma}
\leq M\norm{\fhat}_\infty$ and $\abs{\int_M^s\fhat(\sigma)\,d\sigma}\leq
\norm{\chi_{(M,s)}\fhat}_\infty(s-M)$.  It follows from the
Riemann--Lebesgue lemma that $\Psif(s)=o(s)$ as $\abs{s}\to\infty$.

(d) Note that 
$$
\abs{\Psif(s+h)-\Psif(s)}  =  \left|\intinf \est\left(\frac{1-e^{-iht}}{it}\right)f(t)\,dt\right|
  \leq  C_q\norm{f}_p\abs{h}^{1/p}.
$$
Therefore, $\Psif$ is H\"older continuous on $\R$.
\end{proof}

\begin{remark}
The definition of $\Psif$ depends on the fact that $u_s$ is in the dual space of
$\Lany{p}$, i.e., $u_s\in\Lany{q}$.
However, the dual of the space of Henstock--Kurzweil integrable functions
is the functions of bounded variation, but $u_s$ is not of bounded variation.
In this case, twice integrating the
Fourier transform of a function in $\Lone$ produces a function similar
to $\Psif$ that is of bounded variation.  The
Fourier transform of a Henstock--Kurzweil integrable function can then be
defined as the second distributional derivative.  This case
will be dealt with elsewhere. 
This method does not appear to be applicable to Fourier transforms of functions
in $\Lany{\infty}$, for which a kernel in $\Lone$ is required.
\end{remark}

\begin{remark}
The choice of lower integration limit $0$ in the definition of $\Psif$ is arbitrary.  
If $f\in\Lany{1}$ and $a\in\R$ then
$$
\int_a^s\fhat(\sigma)\,d\sigma -\Psif(s)  =  
\intinf\left(\frac{e^{-iat}-1}{it}\right) f(t)\,dt.
$$
As this is constant, and $\fhat$ is defined by differentiation, 
the choice of lower limit of integration does not affect the definition of the Fourier transform.
With $a$ as the lower limit of integration in the definition of $\Psif$, the
estimate on $\Psif$ then becomes $C_q\norm{f}_p\abs{s-a}^{1/p}$.
\end{remark}

\begin{remark}
If $\mu$ is a finite Borel measure the estimate becomes $\abs{\Psi_{\mu}(s)}\leq\abs{\mu}\abs{s}$
where $\abs{\mu}$ is the total variation of $\mu$.  Indeed, if $\delta_a$ is Dirac measure
supported at $a\in\R$ then $\Psi_{\delta_a}(s)=u_s(a)$.  Then $\widehat{\delta_a}(s)=\Psi'_{\delta_a}(s)=
e^{-ias}$ and this agrees with the distributional definition $\langle\widehat{\delta_a},\phi\rangle=
\langle\delta_a,\phihat\rangle=\phihat(a)$.  Similarly with derivatives of the Dirac distribution.
\end{remark}
\begin{example}
Let $f(x)=\abs{x}^{-1/p}$ where $1<p<\infty$.  A calculation shows 
$\Psif(s)=2\,{\rm sgn}(s)\abs{s}^{1/p}\int_0^\infty
t^{-(1+1/p)}\sin(t)\,dt$.  Hence, $\Psif$ satisfies the inequality in Theorem~\ref{theoremPsif}
but $f$ is not in any $\Lany{p}$ space.  Similarly, if $f$ is taken as an odd function,
except that now the calculation also holds for $p=1$.
\end{example}
\begin{example}
This example shows $\Psif$ need not be of local bounded variation.
Let $f(x)=[J_1(2\sqrt{x})/\sqrt{x}-J_0(2\sqrt{x})]/2$, extended as an odd function to $(-\infty,0)$.
Note that $f\in\Lany{p}$ for $p>4$.
And, $\fhat(s)=i[\sin(1/s)-\cos(1/s)/s]$ \cite[2.13(22), 2.13(25)]{erdelyi}.
Integrating gives $\Psif(s)=is\sin(1/s)$.
Then $\Psif$ is not of bounded variation in any interval containing the origin.
Note that for no $1\leq q\leq\infty$ is $\fhat\in\Lany{q}$.
\end{example}

\section{Fourier transform and the Banach spaces $\Bp$ and $\Ap$}
The Fourier transform will be defined by $\fhat=\Psif'$.  The function
$\Psif$ and its derivative are used to define Banach spaces that
are isometrically isomorphic to $\Lany{p}$.

\begin{theorem}\label{theoremBanachspaces}
Let $1\leq p<\infty$.

(a)  With $\Psif$ given in Theorem~\ref{theoremPsif}
define $\Bp=\{\Psif\mid f\in\Lany{p}\}$.  If $f\in\Lany{p}$ define the norm
$\norm{\Psif}_p'=\norm{f}_p$.   Then $\Bp$ is a Banach space isometrically
isomorphic to $\Lany{p}$.

(b) Define $\Ap=\{\Psif'\mid f\in\Lany{p}\}$.  If $f\in\Lany{p}$ define
$\norm{\Psif'}_p''=\norm{f}_p$.   Then $\Ap$ is a Banach space isometrically
isomorphic to $\Lany{p}$.

(c) The three Banach spaces $\Lany{p}$, $\Bp$ and $\Ap$ are isometrically isomorphic.
\end{theorem}
\begin{proof}
(a) Define the mapping $\Lambda\fn\Lany{p}\to\Bp$ by $\Lambda(f)=\Psif$.  By definition
$\Lambda$ is linear and onto $\Bp$.  To show it is one-to-one suppose $f\in\Lany{p}$
and $\Psif(s)=0$ for all $s\in\R$.  For each
$\phi\in\Sc$ apply the Fubini--Tonelli theorem and integrate by parts,
\begin{align*}
0&=\intinf\Psif(s)\phi'(s)\,ds=\intinf f(t)\intinf\left(\frac{1-e^{-ist}}{it}\right)\phi'(s)
\,ds\,dt\\
 &=-\intinf f(t)\intinf \est\phi(s)\,ds\,dt=-\intinf f(t)\phihat(t)\,dt.
\end{align*}
Note the definition of $u_s$ in the proof of Theorem~\ref{theoremPsif}(a).
The Fubini--Tonelli theorem applies since for $\abs{t}\leq 1$ we can use the estimate
$\abs{u_s(t)}\leq\abs{s}$ from \eqref{u_s} and for $\abs{t}\geq 1$ we can use the estimate
$\abs{u_s(t)}\leq2/\abs{t}$.

Gaussian functions are in the Schwartz space and every Gaussian is the Fourier transform
of a Gaussian.
Therefore,
$\intinf f(t)\Theta_a(x-t)\,dt=0$ for each $a>0$ where $\Theta_a(t)=\frac{e^{-t^2/(4a)}}
{\sqrt{4\pi a}}$ is the heat kernel.  Then $f\ast\Theta_a(x)=0$ for each $x\in\R$ and $a>0$.  But,
$\norm{f-f\ast\Theta_a}_p\to 0$ as $a\to 0^+$.  And, $\norm{f\ast\Theta_a}_p=0$ for each
$a>0$.  Therefore, $f=0$.

(b), (c) The distributional derivative provides a linear isometry from $\Bp$ onto $\Ap$.
It is one-to-one since, due to the estimate $\abs{\Psif(s)}\leq C_q\norm{f}_p\abs{s}^{1/p}$ from
Theorem~\ref{theoremPsif}(a), the only constant in $\Bp$ is $0$.  Hence, the spaces
$\Ap$ and $\Bp$ are isometrically isomorphic.
\end{proof}

Now we can define the Fourier transform of functions in $\Lany{p}$ for $1\leq p<\infty$.

\begin{defn}\label{defnfhatLp}
Let $1\leq p<\infty$.  Let $f\in\Lany{p}$. 
Define the Fourier transform as $\FT\fn\Lany{p}\to\Ap$ by $\fhat=\Psif'$,
i.e., 
$\langle \fhat,\phi\rangle=-\langle \Psif,\phi'\rangle=-\intinf\Psif(s)\phi'(s)\,ds$ 
for each $\phi\in\Sc$.
\end{defn}
The space of Fourier transforms on $\Lany{p}$ is then denoted $\Ap$.

\section{Properties of the Fourier transform}
The Fourier transform
of a tempered distribution $T$ is $\langle \hatT,\phi\rangle=\langle T, \phihat\rangle$
for each Schwartz function $\phi$.
Since $\Psif$ is continuous and of slow (polynomial) growth then $\fhat$ as defined
by Definition~\ref{defnfhatLp} is a tempered distribution.
Hence, Fourier transforms in 
$\Ap$ inherit the properties of Fourier transforms in $\Scp$.

\begin{theorem}\label{theoremequivalentdefns}
Definition \ref{defnfhatLp} agrees with the definition for the Fourier transform of a tempered
distribution.
\end{theorem}
\begin{proof}
This is contained in the proof of Theorem~\ref{theoremBanachspaces}(a).
\end{proof}

Using Definition~\ref{defnfhatLp}
many of the usual properties of Fourier transforms have a different form,
cf. \cite[\S30]{donoghue} or \cite[2.3.22]{grafakos}.
The same notation is used below.

\begin{theorem}\label{theoremalgebraic}
Let $1\leq p<\infty$.  Let $f\in\Lany{p}$ and define $\fhat$ with 
Definition~\ref{defnfhatLp}.
\begin{enumerate}
\item[(a)]
The Fourier transform is a bounded linear isomorphism and isometry, $\FT\fn\Lany{p}\to\Ap$,
$\norm{\fhat}''_p=\norm{f}_p$.
\item[(b)]
Let $F^{(n-1)}$ be absolutely continuous such that $F^{(n)}\in\Lany{p}$.
For each $0\leq k\leq n-2$ suppose
$\llim_{\abs{t}\to\infty}F^{(k)}(t)/t=0$.
Then $\widehat{D^nF}=\Psi'_{\!F^{(n)}}$ where
$$
\Psi_{\!F^{(n)}}(s)=\intinf\frac{n!}{it^{n+1}}\left[1-e^{-ist}\sum_{k=0}^n\frac{(ist)^k}{k!}\right]
F(t)\,dt.
$$
\item[(c)]
$D^n\fhat=\Psif^{(n+1)}$.
\item[(d)]
Let $a\in\R$.  Then $\widehat{\tau_af}=\Psi_{\tau_af}'$ where
$\Psi_{\tau_af}(s)=f\ast\tilde{u}_s(-a)$.
\item[(e)]
Let $F(t)=e^{iat}f(t)$.  Then
$\hat{F}=(\tau_a\Psif)'$.
\item[(f)]
$\hat{{\tilde f}}=\Psi_{\!\tilde{f}}'$ where $\Psi_{\!\tilde{f}}=-\tilde{\Psi}_f$.
\item[(g)]
If $\fhat\in\Lany{r}$ for some $1\leq r<\infty$ then $\hat{\fhat}=\Psi_{\!\fhat}'$.
\item[(h)]
Let $g(x)=f(ax+b)$ for $a,b\in\R$ with $a\neq 0$.  Then
$\ghat=\Psi_{\!g}'$ where $\Psi_{\!g}(s)={\rm sgn(a)}\Psi_{\!\tau_{-b}f}(s/a)$.
\end{enumerate}
\end{theorem}

\begin{proof}
\begin{enumerate}
\item[(a)]
This follows from Definition~\ref{defnfhatLp} and Theorem~\ref{theoremBanachspaces}.
\item[(b)]
We have $F^{(n-1)}(x)=F^{(n-1)}(0)+\int_0^xF^{(n)}(t)\,dt$ 
so the H\"older inequality gives
$\abs{F^{(n-1)}(x)}\leq\abs{F^{(n-1)}(0)}+\norm{F^{(n)}}_p\abs{x}^{1/q}$.  It follows that
$\lim_{\abs{t}\to\infty}F^{(n-1)}(t)/t=0$.  Induction on $n$ establishes the formula
$$
u_s^{(n)}(t)=\frac{(-1)^nn!}{it^{n+1}}\left[1-e^{-ist}\sum_{k=0}^n\frac{(ist)^k}{k!}\right].
$$
Notice that $u_s$ and its derivatives are continuous and of order $1/t$ as $\abs{t}\to\infty$.
An $n$-fold integration by parts can then be performed on $\Psi_{\!F^{(n)}}(s)=\intinf
u_s(t)F^{(n)}(t)\,dt$.
\item[(c)]
$\langle D^n\fhat,\phi\rangle=\langle\Psif^{(n+1)},\phi\rangle=(-1)^{n+1}\langle\Psif,\phi^{(n+1)}
\rangle$.  Integration by parts gives another formula for the derivative of a Fourier transform.
\item[(d)]
Follows immediately from the definition.
\item[(e)]
We have
$$
\Psi_{\!F}(s)=\intinf\left(\frac{1-e^{-i(s-a)t}}{it}-\frac{1-e^{iat}}{it}\right)f(t)\,dt=
\Psif(s-a)-\Psif(-a)
$$
so that $\hat{F}=(\tau_a\Psif)'$.
\item[(f)]
We have
$$
\Psi_{\!\tilde{f}}(s)=\intinf\left(\frac{1-e^{-ist}}{it}\right) f(-t)\,dt
=-\intinf\left(\frac{1-e^{ist}}{it}\right) f(t)\,dt.
$$
\item[(g)]
Follows immediately from the definition.
\item[(h)]
Change variables in the integral defining $\Psi_{\!g}$.
\end{enumerate}
\end{proof}
\begin{remark}
For use in (g),
there are various conditions that ensure $\fhat\in\Lany{r}$.
Three such conditions are as follows.

(a) Suppose $f\in\Lany{p}$ for some $1< p\leq 2$.  By the
Hausdorff--Young--Babenko--Beckner inequality, $\fhat\in\Lany{q}$ where $2\leq q<\infty$
is the conjugate exponent.

(b) Suppose $f$ is absolutely continuous and $f,f'\in\Lone$.  This implies $\lim_{\abs{x}\to\infty}f(x)=0$
\cite[Lemma~2]{talvilafouriermaa}.  
And, $\fhat$ and $\widehat{f'}$ are continuous and bounded.
Integration by parts shows $\fhat(s)=\widehat{f'}(s)/(is)$.  Hence, $\fhat\in\Lany{r}$
for each $1<r<\infty$.

(c) Suppose $f\in\Lone$ and $f$ is of bounded variation such that $\lim_{\abs{x}\to\infty}
f(x)=0$.  Integrate by parts to get $\fhat(s)=(1/(is))\intinf e^{-ist}\,df(t)$.
Then $\abs{\fhat(s)}\leq V\!f/\abs{s}$.  Since $\fhat$ is continuous, $\fhat\in\Lany{r}$
for each $1<r<\infty$.
\end{remark}

\section{Integration in $\Ap$}
If $I$ is a closed interval in $\R$,
denote the functions of bounded variation on $I$ by $\bv(I)$.
If $g\in\bv(\R)$ then $g(\infty)$ is defined as $\lim_{x\to\infty}g(x)$.
Similarly for $g(-\infty)$.

Since $\Psif$ is continuous $\fhat$ is integrable in the sense of the
continuous primitive integral.
In this theory, if $f=F'$ for
a function $F$ that is continuous on $[a,b]$ 
then $\intab F'g=F(b)g(b)-F(a)g(a)-\intab F(x)\,dg(x)$ where $g\in\bv([a,b])$
and the last integral is a Henstock--Stieltjes integral.  
See \cite[Definition~6]{talviladenjoy}.
Following notation there, Lebesgue integrals will explicitly display the variable
and measure, continuous primitive integrals will show neither.

The growth estimates for $\Psif$ in Theorem~\ref{theoremPsif} lead to the following
definitions of integrals for distributions in $\Ap$.

\begin{defn}\label{defnintegrationinAp}
Let $f\in\Lany{p}$ for some $1\leq p<\infty$.

(a) Let $-\infty<a<b<\infty$.
Let $g\in\bv([a,b])$.
Then $\intab\fhat g=\intab\Psif'g=\Psif(b)g(b)-\Psif(a)g(a)-\intab\Psif(t)\,dg(t)$.

(b) Let $g\in\bv(\R)$ such that $g(x)=o(\abs{x})^{-1/p}$ as $\abs{x}\to\infty$.
If $p=1$ then $g(x)=O(x^{-1})$ as $\abs{x}\to\infty$ suffices.
Then $\intinf\fhat g=\intinf\Psif'g=-\intinf\Psif(t)\,dg(t)$.

(c) Let $g\in\bv([\delta,a])$ for each $0<\delta<a<\infty$ such that $g(x)=o(x^{-1/p})$ as $x\to0^+$.  
Then
$\int_0^a\fhat g=\int_0^a\Psif'g=\Psif(a)g(a)-\int_0^a\Psif(t)\,dg(t)$.
In particular, $\int_a^b\fhat=\Psif(b)-\Psif(a)$ for all $a,b\in\R$.
\end{defn}

\begin{remark}\label{remarkBabenko}
Definition~\ref{defnfhatLp} and Theorem~\ref{theoremPsif} give an estimate on the integrals of $\fhat$,
$\abs{\int_0^s\fhat}\leq C_q\norm{f}_p\abs{s}^{1/p}$.  If $1\leq p\leq 2$ then the
Hausdorff--Young--Babenko--Beckner inequality is
$\norm{\fhat}_q\leq B_q\norm{f}_p$ where the sharp constant is
$B_q=(2\pi)^{1/q}[q^{1-2/q}(q-1)^{1/q-1}]^{1/2}=(2\pi)^{1/q}(p^{1/p}/q^{1/q})^{1/2}$.
See \cite{babenko}, \cite{beckner} and \cite{liebloss}.  By the H\"older inequality,
$\abs{\int_0^s\fhat(\sigma)\,d\sigma}\leq \norm{\fhat}_q\abs{s}^{1/p}\leq
B_q\norm{f}_p\abs{s}^{1/p}$.  It is then interesting to compare the values of $B_q$
and $C_q$ in these inequalities.  Note that $B_\infty=C_\infty=1$ and $B_2=C_2=\sqrt{2\pi}$.
As well, $\lim_{q\to1^+}C_q=\infty$ while $B_1=2\pi$.
The values of $B_q$ can be computed explicitly but $C_q$ is known exactly only for
$q=2,4,6,\ldots$.  Due to numerical calculations, we conjecture that $C_q>B_q$ for $1<q<2$
while $C_q<B_q$ for $2<q<\infty$.  Hence, for $1<p<2$, if the conjecture holds then the inequality 
$\abs{\int_0^s\fhat(\sigma)\,d\sigma}\leq C_q\norm{f}_p\abs{s}^{1/p}$ is sharper than that obtained
in this manner using the Hausdorff--Young--Babenko--Beckner inequality.
\end{remark}

\section{Exchange formula and inversion}\label{sectionexchange}
The exchange formula is $\intinf \fhat(x)g(x)\,dx=\intinf f(x)\hat{g}(x)\,dx$, which
holds for $f,g\in\Lany{p}$ ($1\leq p\leq 2$).  
If $T$ is a tempered distribution and $\phi$ is a Schwartz function
then the formula is just the definition of the Fourier transform,
$\langle \hatT,\phi\rangle=\langle T,\phihat\rangle$.  This can now
be extended to $\Ap$.

\begin{theorem}\label{theoremexchange}
Let $1\leq p<\infty$ with $q$ the conjugate exponent. Let $f\in\Lany{p}$.
Suppose $g\fn\R\to\R$ is of bounded variation with limit $0$ at $\pm\infty$
such that $\intinf\abs{s}^{1/p}\abs{dg(s)}<\infty$.  Then $\ghat\in\Lany{q}$ and
$\intinf\fhat g=\intinf f(s)\ghat(s)\,ds$.
And, $\abs{\intinf\fhat g}\leq C_q\norm{f}_p\intinf\abs{s}^{1/p}\abs{dg(s)}$,
with $C_q$ as in Theorem~\ref{theoremPsif}.
\end{theorem}
\begin{proof}
It suffices to write $g=g_1-g_2$ where
$g_1,g_2\fn[0,\infty]\to\R$ and decrease to $0$.
Then $\int_0^\infty\abs{s}^{1/p}\abs{dg(s)}=-\int_0^\infty s^{1/p}\,dg_1(s)-\int_0^\infty s^{1/p}\,dg_2(s)$.
Note that $\lim_{x\to\infty}\int_x^\infty s^{1/p}\,dg_1(s)=0$.
For $x>0$, $-\int_x^\infty s^{1/p}\,dg_1(s)\geq x^{1/p}g_1(x)\geq 0$.  Hence, $g(s)=o(\abs{s}^{-1/p})$ as
$\abs{s}\to\infty$.

Integrate by parts,
\begin{eqnarray*}
\intinf\fhat g  & = &  \intinf\Psif'g =-\intinf\Psif(s)\,dg(s)=-\intinf\intinf \left(\frac{1-e^{-ist}}{it}\right) f(t)\,dt\,dg(s)\\
 & = & -\intinf f(t)\intinf\left(\frac{1-e^{-ist}}{it}\right)dg(s)\,dt.
\end{eqnarray*}
The Fubini--Tonelli theorem applies since $\intinf \abs{\Psif(s)}\abs{dg(s)}\leq
C_q\norm{f}_p\intinf\abs{s}^{1/p}\abs{dg(s)}<\infty$.
Then, upon integrating by parts again,
\begin{eqnarray*}
\intinf\fhat g  & = & -\intinf f(t)\left\{\left[\left(\frac{1-e^{-ist}}{it}\right)g(s)\right]_{s=-\infty}
^\infty-\intinf e^{-ist}g(s)\,ds\right\}dt\\
 & = & \intinf f(t)\ghat(t)\,dt.
\end{eqnarray*}
The estimate
$$
\left|\intinf f(t)\ghat(t)\,dt\right|=\left|\intinf\fhat g\right|\leq
C_q\norm{f}_p\intinf\abs{s}^{1/p}\abs{dg(s)}<\infty
$$
holds for each $f\in\Lany{p}$ so $\ghat\in\Lany{q}$.
\end{proof}

\begin{example}\label{exampleexchangesin(x)/x}
Let $g(x)=\sin(x)/x$.  Then $\ghat(s)=\pi\chi_{(-1,1)}(s)$ so $\ghat\in\Lany{q}$ for each $1\leq q\leq\infty$.
But, $g$ is not of bounded variation.  And,
$\intinf\abs{s}^{1/p}\abs{dg(s)}=2\int_0^\infty s^{1/p}\abs{g'(s)}\,ds=2\int_0^\infty s^{1/p}\abs{
\cos(s)/s-\sin(s)/s^2}\,ds =\infty$.  Hence, none of the conditions on $g$ in Theorem~\ref{theoremexchange}
hold.  Those conditions are then sufficient but not necessary.
\end{example}

\begin{example}
Let $1<p<\infty$.
Let $g(x)=x^{-\alpha}$ for $x>1$ and $g=0$, otherwise.  Here, $1/p<\alpha<1$.
Then $g\in\Lany{p}\cap\bv(\R)$.  And, $\intinf\abs{s}^{1/p}\abs{dg(s)}=\alpha\int_1^\infty s^{-(\alpha+1-1/p)}\,ds<\infty$.
For $s>0$ note that $\ghat(s)=s^{\alpha-1}\int_s^\infty e^{-it}t^{-\alpha}\,dt\sim s^{\alpha-1}\int_0^\infty
e^{-it}t^{-\alpha}\,dt$ as $s\to0^+$.  Integration by parts shows
$\ghat(s)=e^{-is}/(is)+O(s^{-2})$ as $s\to\infty$.  It follows that $\ghat\in\Lany{q}$.  Note that
$g\not\in\Lone$.
\end{example}

The exchange formula leads to inversion in the $\Lany{p}$ norm using a summability kernel.

\begin{theorem}\label{theoreminversion}
Let $1\leq p<\infty$.  Let $f\in\Lany{p}$.
Let $\psi\in\Lone$ such that $\intinf\psi(t)\,dt=1$,
$\psihat$ is absolutely continuous,
$\intinf\abs{s}^{1/p}\abs{\psihat(s)}\,ds<\infty$ and $\intinf\abs{s}^{1/p}\abs{\psihat\,'(s)}\,ds<\infty$.
For each $a>0$ define $\psi_a(x)=\psi(x/a)/a$.  Define a family of kernels $K_a(s)=\widehat{\psi_a}(s)$.
Let $e_x(s)=e^{ixs}$ and $I_a[f](x)=(1/(2\pi))\intinf e_{x}K_a\fhat$.
Then $\lim_{a\to0^+}\norm{f-I_a[f]}_p=0$.
\end{theorem}
\begin{proof}
The proof of Theorem~\ref{theoremexchange} shows that theorem applies to the function
$s\mapsto  e^{ixs}K_a(s)$.  Then $I_a[f](x)=(1/(2\pi))\intinf f(t)\widehat{K_a}(t-x)\,dt
=f\ast\psi_a(x)$.
By approximation with convolutions in norm 
\cite[Theorem~8.14]{folland} the conclusion follows.
\end{proof}

\begin{example}
The three commonly used summability kernels:
Ces\`aro--Fej\'er ($K_a(s)=(1-a\abs{s})\chi_{[-1/a,1/a]}(s)/(2\pi)$ with 
${\widehat K_a}(t)=2a\sin^2[t/(2a)]/(\pi t^2)$),
Abel--Poisson ($K_a(s)=e^{-a\abs{s}}/\pi$ with ${\widehat K_a}(t)=2a/[\pi(t^2+a^2)]$) and
Gauss--Weierstrass ($K_a(s)=e^{-a^2s^2}/(2\pi)$ with ${\widehat K_a}(t)=e^{-t^2/(4a^2)}/(2\sqrt{\pi}a)$)
satisfy the conditions of Theorem~\ref{theoreminversion}.
While the Dirichlet kernel ($K_a(s)=\chi_{[-1/a,1/a]}(s)/(2\pi)$ with ${\widehat K_a}(t)=\sin(t/a)/(\pi t)$)
does not.

The way summability kernels are defined in Theorem~\ref{theoreminversion} means
the inverses of these kernels are associated with the given names.
\end{example}

\section{Convolution}\label{sectionconvolution}
The exchange formula is used with Fourier transforms of convolutions.
\begin{theorem}\label{theoremconvolution}
Let $f\in\Lany{p}$ for some $1\leq p<\infty$.
Let $g_1\in\Lone$ such that $\gonehat$ is of bounded variation.
Let $g_2\in\Lone$ such that $g_2$ is of bounded variation with limit $0$ at $\pm\infty$
and $\intinf\abs{s}^{1/p}\abs{dg_2(s)}<\infty$.
Then $f\ast g_1\in\Lany{p}$ and $\intinf\widehat{f\ast g_1}g_2=\intinf\fhat\gonehat
g_2$.
\end{theorem}
\begin{proof}
By Young's inequality \cite[8.7]{folland},
$f\ast g_1\in\Lany{p}$ so Theorem~\ref{theoremexchange} gives
$\intinf \widehat{f\ast g_1}g_2=\intinf f\ast g_1(t)\gtwohat(t)\,dt$.
Note that 
\begin{equation}
\intinf\abs{s}^{1/p}\abs{d(\gonehat g_2(s))}
\leq \intinf\abs{s}^{1/p}\abs{\gonehat(s)}\abs{dg_2(s)}
+\intinf\abs{s}^{1/p}\abs{g_2(s)}\abs{d\gonehat}.\label{gonehat}
\end{equation}
The first integral on the right side of \eqref{gonehat} is bounded since
$\gonehat$ is bounded.  The proof of Theorem~\ref{theoremexchange} shows
$\abs{s}^{1/p}g_2(s)$ is bounded.  Hence, the left side of \eqref{gonehat}
is bounded.  Theorem~\ref{theoremexchange} then shows that
$\intinf\fhat\gonehat
g_2=\intinf f(t)\widehat{\gonehat g_2}(t)\,dt$.
And,
\begin{eqnarray*}
\widehat{\gonehat g_2}(t) & = & \intinf e^{-ist}\intinf e^{-isu}g_1(u)g_2(s)\,du\,ds
 =  \intinf g_1(u)\intinf e^{-i(u+t)s}g_2(s)\,ds\,du\\
 & = & \intinf g_1(u-t)\widehat{g_2}(u)\,du.
\end{eqnarray*}
Since $g_1,g_2\in\Lone$ the Fubini--Tonneli theorem permits interchange of integrals.

Now, 
\begin{eqnarray*}
\intinf f(t)\widehat{\gonehat g_2}(t)\,dt & = & \intinf f(t)\intinf g_1(u-t)\widehat{g_2}(u)\,du\,dt\\
 & = & \intinf\intinf f(t)g_1(u-t)\,dt\,\widehat{g_2}(u)\,du\\
 & = & \intinf f\ast g_1(u)\widehat{g_2}(u)\,du.
\end{eqnarray*}
By Theorem~\ref{theoremexchange}, $\widehat{g_2}\in\Lany{q}$ with $q$ the conjugate
of $p$ so the Fubini--Tonneli theorem applies.
\end{proof}

Theorem~\ref{theoremconvolution} requires $\gonehat$ be of bounded variation.
The following proposition provides sufficient conditions for a Fourier transform
to be in $\Lone$.  Further conditions then show when a Fourier transform is of
bounded variation.

\begin{prop}\label{propvariation}
(a) Let $g\in\Lone$ be absolutely continuous such that $\lim_{\abs{x}\to\infty}g(x)=0$
and $g'\in\Lany{p}$ for some $1<p\leq 2$.  Then $\ghat\in\Lone$.

(b) Let $g\in\Lone$.  Define $h(x)=xg(x)$.  Suppose $h\in\Lone$ such that
$h$ is absolutely continuous, $\lim_{\abs{x}\to\infty}h(x)=0$ and $h'\in\Lany{p}$
for some $1<p\leq 2$.  Then $\ghat\in\bv(\R)$.

(c) Let $g$ be absolutely continuous such that $g,g'\in\Lone$ and $g'\in\bv(\R)$.
Then $\ghat\in\Lone$.

(d) Let $g\in\Lone$.  Define $h(x)=xg(x)$.  Suppose $h\in\Lone$ such that
$h$ is absolutely continuous, $h'\in\Lone$ and $h'\in\bv(\R)$.
Then $\ghat\in\bv(\R)$.
\end{prop}
\begin{proof}
(a) Since $\ghat$ is continuous it is locally integrable.  For $s\not=0$, integrate
by parts to get $\ghat(s)=\widehat{g'}(s)/(is)$.  Then, using the 
H\"older and
Hausdorff--Young--Babenko--Beckner inequalities,
$$
\int_1^\infty\abs{\ghat(s)}\,ds  =  \int_1^\infty\abs{\widehat{g'}(s)}\frac{ds}{s}
  \leq  \norm{\widehat{g'}}_q\left(\int_1^\infty s^{-p}\,ds\right)^{1/p}
  \leq  \frac{\sqrt{2\pi}\norm{g'}_p}{(p-1)^{1/p}}.
$$
cf. \cite[p. 262, Exercise~24]{folland}

(b) The proof as per (a).

(c) The conditions $g$ absolutely continuous with $g,g'\in\Lone$ imply that
$\lim_{\abs{x}\to\infty}g(x)=0$ \cite[Lemma~2]{talvilafouriermaa}.
The conditions $g'\in\Lone$ and of bounded variation imply that
$\lim_{\abs{x}\to\infty}g'(x)=0$.  Integration by parts then shows
$\ghat(s)=-(1/s^2)\intinf e^{-ist}\,dg'(t)$, from which
$\abs{\ghat(s)}\leq V\!g'/s^2$ for $s\not=0$.  Since $\ghat$ is continuous the result
follows.

(d) The proof as per (c).
\end{proof}

\begin{remark}\label{remarkbv}
The conditions for $\ghat$ being of bounded variation are not necessary.
For example, $g(x)=\sin(x)/x$ does not satisfy the conditions in (b) or (d)
but $\ghat=\pi\chi_{(-1,1)}$ which is of bounded variation.

If  $1/2\leq\alpha<1$ and
$$
g(x)=\left\{
\begin{array}{cl}
x^{-\alpha}, & 0<x\leq 1\\
e^{1-x}, & x\geq 1\\
0, & x\leq 0,
\end{array}
\right.
$$
then $g$ satisfies the conditions in (b) but not those of (d).
\end{remark}

\begin{theorem}\label{theoremg1g2convolution}
Let $f\in\Lany{p}$ for some $1\leq p<\infty$.
Suppose $g_1\in\Lone$ is of bounded variation 
such that $\intinf\abs{s}^{1/p}\abs{dg_1(s)}<\infty$.
Suppose $g_2\in\Lone$ such that $\intinf\abs{s}^{1/p}\abs{g_2(s)}\,ds<\infty$.
Then
$\intinf\fhat g_1\ast g_2=\intinf f(s)\gonehat(s)\gtwohat(s)\,ds$.
\end{theorem}
\begin{proof}
Show that $g_1\ast g_2$ satisfies the conditions of Theorem~\ref{theoremexchange}.
If $N\in\N$ and $-\infty=s_0<s_1<\ldots<s_N=\infty$ then
$$
\sum_{n=1}^N\abs{g_1\ast g_2(s_{n})-g_1\ast g_2(s_{n-1})}
\leq \intinf\sum_{n=1}^N\abs{g_1(s_{n}-t)-g_1(s_{n-1}-t)}\abs{
g_2(t)}\,dt.
$$
And, $V\!g_1\ast g_2\leq V\!g_1\norm{g_2}_1$.

Since $g_1\in\Lone$ is of bounded variation it has limit $0$ at $\pm\infty$.
Dominated convergence shows $g_1\ast g_2$ has limit $0$ at $\pm\infty$.

Now suppose $z_n\in[s_{n-1},s_{n}]$.  Consider the Riemann sums,
\begin{align}
&\sum_{n=1}^N\abs{z_n}^{1/p}\abs{g_1\ast g_2(s_{n})-g_1\ast g_2(s_{n-1})}\notag\\
& \leq  \sum_{n=1}^N\intinf\abs{(z_n-t)+t}^{1/p}\abs{g_1(s_{n}-t)-g_1(s_{n-1}-t)}\abs{
g_2(t)}\,dt\notag\\
& \leq  \intinf\sum_{n=1}^N\abs{z_n-t}^{1/p}\abs{g_1(s_{n}-t)-g_1(s_{n-1}-t)}\abs{
g_2(t)}\,dt\notag\\
&\qquad+ \intinf\abs{t}^{1/p}\sum_{n=1}^N\abs{g_1(s_{n}-t)-g_1(s_{n-1}-t)}\abs{
g_2(t)}\,dt\notag\\
&\leq \sup_{t\in\R}\sum_{n=1}^N\abs{z_n-t}^{1/p}\abs{g_1(s_{n}-t)-g_1(s_{n-1}-t)}
\norm{g_2}_1+V\!g_1\intinf\abs{t}^{1/p}\abs{
g_2(t)}\,dt.\label{g1g2convolution}
\end{align}
For each $t\in\R$, the
tagged partition can be chosen so that the Riemann sum in \eqref{g1g2convolution}
approximates the integral $\intinf\abs{s-t}^{1/p}\abs{dg_1(s-t)}
=\intinf\abs{s}^{1/p}\abs{dg_1(s)}$ within any given
$\epsilon>0$.  See, for example, \cite{mcleod}.
This shows that $\intinf\abs{s}^{1/p}\abs{d(g_1\ast g_2)(s)}$ exists.

Use of Theorem~\ref{theoremexchange} completes the proof.
\end{proof}

\end{document}